\documentclass[12pt,leqno]{article}

\usepackage{amssymb,amsmath,amsthm}

\newtheorem{theorem}{Theorem}[section]
\newtheorem{lemma}[theorem]{Lemma}

\newtheorem{corollary}[theorem]{Corollary}

\theoremstyle{definition}

\theoremstyle{remark}
\newtheorem{remark}[theorem]{Remark}

\numberwithin{equation}{section}

\newcommand{\wedgeco}%
{\displaystyle\operatornamewithlimits{\wedge}_{\raise3pt\hbox{,}}}

\newcommand{\ostar}{\odot\kern-6.4pt\ast}

\def\zs#1{_{\lower2pt\hbox{$\scriptstyle#1$}}}

\newcommand{\hb}{\hbar}

\newcommand{\pa}{\partial}
\newcommand{\ve}{\varepsilon}

\newcommand{\wh}{\widehat}
\newcommand{\wt}{\widetilde}
\newcommand{\od}{\overset{\text{\rm def}}{=}}

\newcommand{\Exp}{\operatorname{Exp}}
\newcommand{\Graph}{\operatorname{Graph}}

\newcommand{\id}{\operatorname{id}}

\newcommand{\bR}{\mathbb{R}}

\newcommand{\cE}{\mathcal{E}} 
 
\newcommand{\cH}{\mathcal{H}}

\newcommand{\cX}{\mathfrak{X}} 
 
\newcommand{\cS}{\mathcal{S}}

\begin{document}

\title{Intrinsic Dynamics of Symplectic Manifolds: 
Membrane Representation and Phase Product}  

\author{Mikhail Karasev\thanks{This research was supported 
in part by the Russian Foundation for Basic Research, 
Grant~02-01-00952.}\\
\\
\small Moscow Institute of Electronics and Mathematics\\
\small Moscow 109028, Russia\\
\small karasev@miem.edu.ru}

\date{}

\maketitle

\begin{abstract}
On general symplectic manifolds 
we describe a correspondence between symplectic transformations 
and their phase functions. 
On the quantum level, this is a correspondence between unitary
operators and phase functions of the WKB-approximation. 
We represent generic functions via symplectic area of membranes and
consider related geometric properties of the noncommutative
phase product. 
An interpretation of the phase product 
in terms of symplectic groupoids 
and the groupoid extension of Lagrangian submanifolds
are described. 
The membrane representations 
of corresponding Lagrangian phase functions are obtained.  
This paper uses the intrinsic dynamic approach 
based on the notion of Ether Hamiltonian which is a
generalization of the notion of symplectic connection. 
We demonstrate that this approach works for torsion case as well.
\end{abstract}

\maketitle

\section{Introduction}

The {\it intrinsic dynamics\/} of symplectic manifolds is
generated by symplectic connections, 
but this dynamics is not equivalent to the kinematics  
of parallel translations of vectors or tensors by means
of the connection.  
The intrinsic dynamics is based on a certain object 
(called the Ether Hamiltonian~\cite{K})
which can be considered as a connection in the function
bundle over the manifold. 
This object  allows one to translate functions, rather than
vectors, and thus to exploit certain hidden geometry
which we never see in the habitual analysis produced 
by the infinitesimal geometrical tools.

In the given paper, we apply these ideas to investigate 
the following question:
how to associate functions on a phase space with symplectic
transformations of this space and what is the product of
functions corresponding to the composition of transformations? 
On the Euclidean phase spaces with the canonical structure 
$dp\wedge dq$ this is a routine question about phase functions
(or generating functions) of symplectic
transformations~\cite{A,H}.  
The choice of a phase function depends on the choice of
polarization on the double phase space. 

In the general case, the given Ether Hamiltonian
on the symplectic manifold determines 
a natural choice of polarization, 
and thus one obtains a correspondence between functions and
symplectic transformations. 

Moreover, each function close enough to constant
can be represented as the symplectic area of a membrane whose
boundary is organized by means of Ether
geodesics and the related symplectic transformation.
Composition of transformations makes up the composition of
membranes. Via the Stokes theorem, the area of this composition
represents a noncommutative product of functions.

The corresponding composition of Lagrangian submanifolds 
is generated by the groupoid multiplication operation. 
This operation can also be used for transformations and extensions of
Lagrangian submanifolds, and for various types 
of membrane representations of generating functions 
related to these submanifolds.

In the last section we consider a generalization
of the results to the torsion case.

\section{Ether Hamiltonian}

Let $\cX$ be a manifold with symplectic form $\omega$ 
and symplectic torsion free connection~$\Gamma$ 
(that is, $\omega$ is covariantly constant in the sense
of~$\Gamma$).  

The connection generates the kinematic geometry on~$\cX$: 
parallel translations of vectors, geodesics, exponential
mappings, etc. 
There are no forces in this geometry, no changing of momenta, 
and no opportunities to translate functions or curves and
surfaces. 
That is why one needs to integrate the pair $(\omega,\Gamma)$ up
to a more substantial object which we call the {\it Ether
Hamiltonian}. 
It generates a dynamic geometry on~$\cX$. 
We recall some notions from~\cite{K}.

The Ether Hamiltonian is a $1$-form of~$\cX$ with values 
in $C^\infty(\cX)$ satisfying the zero curvature equation 
and some boundary and skew-symmetry conditions.
We use the following notation for this Hamiltonian:
$$
\cH_x(z)=\sum^{2n}_{j=1}\cH_x(z)_j dx^j,\qquad
2n=\dim \cX,\quad x,z\in\cX.
$$
The {\it zero curvature equation\/} is 
\begin{equation}
\pa\cH+\frac12\{\cH\wedgeco\cH\}=0.
\tag{2.1}
\end{equation}
Here the Poisson brackets $\{\cdot,\cdot\}$ are taken with
respect to the symplectic form $\omega=\omega(z)$, 
and $\pa =\pa_x$ denotes the differential of a form  
at the point~$x$. 

The {\it boundary conditions\/} are 
\begin{equation}
\cH_x(z)\Big|_{x=z}=0,\quad 
D\cH_x(z)\Big|_{x=z}=2\omega(z),\quad
D^2\cH_x(z)\Big|_{x=z}=2\omega(z)\Gamma(z).
\tag{2.2}
\end{equation}
Here $D=D_z$ is the derivative with respect to~$z$.
Note that everywhere we use identical notations 
both for a differential $2$-form itself 
and for the matrix of coefficients of this form 
in local coordinates. 
In (2.2), of course, $\omega(z)$ stays for the matrix of
coefficients of the form~$\omega$ at the point~$z$, 
and $\Gamma(z)$ is the matrix of Christoffel symbols 
of the connection~$\Gamma$ at the point~$z$.

Let $z\to s_x(z)$ be a trajectory of
the Ether Hamiltonian  
(that is, the ``time'' derivative $\pa s_x(z)$ coincides
with the Hamiltonian vector field of $\cH_x$ 
at the point $s_x(z)$), 
and the ``initial'' data are $s_x(z)\big|_{x=z}=z$,
see (3.4) below.
  
The {\it skew-symmetry condition\/} is 
\begin{equation}
\cH_x\big(s_x(z)\big)=-\cH_x(z).
\tag{2.3}
\end{equation}

For each $x\in\cX$ and $v\in T_x\cX$, 
by $\Exp_x(vt)$ we denote the Hamilton trajectory on~$\cX$ which
corresponds to the Hamilton function $\frac12v\cH_x$ and starts
at~$x$ when $t=0$.

\begin{theorem}
{\rm(i)} In a neighborhood of the diagonal $z=x$,
there exists a solution of Eq.~{\rm(2.1)} 
satisfying conditions {\rm(2.2), (2.3)}.

{\rm(ii)} Mappings $s_x$ are symplectic transformations
of~$\cX$. 
The point~$x$ is an isolated fixed point of~$s_x${\rm:}
\begin{equation}
s_x(x)=x.
\tag{2.4}
\end{equation}
The mapping $s_x$ is an involution{\rm:}
\begin{equation}
s^2_x=\id.
\tag{2.5}
\end{equation}
Moreover, the family $\{s_x\}$ is related to the connection $\Gamma$
by the formula 
\begin{equation}
\Gamma(z)=-\frac12 D^2 s_x(z)\Big|_{x=z}.
\tag{2.6}
\end{equation}

{\rm(iii)} If $\{s_x(z)\mid x\in \cX\}$ is a smooth family of
symplectic transformations of~$\cX$ satisfying {\rm(2.4),
(2.5)}, then formula {\rm(2.6)} determines a symplectic
connection on~$\cX$, and formula 
\begin{equation}
\cH_x(z)=\int^z_x\langle\pa s_x(s_x(z)),\omega(z)\,dz \rangle
\tag{2.7}
\end{equation}
determines the Ether Hamiltonian {\rm(}the solution of
{\rm(2.1)--(2.3))}. 

{\rm(iv)}
The mappings $\Exp_x$ are related to~$\cH_x$ and 
to the family $\{s_x\}$ as follows{\rm:}
\begin{equation}
\cH_x(\Exp_x(v))=-\cH_x(\Exp_x(-v)),\qquad
s_x(\Exp_x(v))=\Exp_x(-v).
\tag{2.8}
\end{equation}
\end{theorem}

We call the mapping  $s_x$ a {\it reflection}, 
and we call $\Exp_x$ an Ether {\it exponential mapping}. 
The curve  
$\{\Exp_x(vt)\mid -\ve<t<\ve\}$ is called 
the {\it Ether geodesics through the mid-point}~$x$.

The composition of two reflections 
$g_{x,y}=s_x\circ s_y$ 
we call the {\it Ether translation}. 
This map coincides with the shift along trajectories of the
Ether dynamical system with Hamiltonian $\cH_x$ 
when the ``time'' varies from~$y$ to~$x$ (see in~\cite{K}).

In general, the reflections~$s_x$ do not preserve 
the connection~$\Gamma$, 
the mappings $\Exp_x$ do not coincide with the exponential
mappings $\exp_x$ generated by~$\Gamma$,
and so the Ether geodesics do not coincide with the 
$\Gamma$-geodesics.

\section{Phase functions}

In the double $\cX\times\cX$ with the symplectic structure 
$\omega\ominus\omega$ we have the Lagrangian fibration
$\cS=\{\cS_x\mid x\in\cX\}$ whose fibers are graphs of
reflections: $\cS_x=\Graph(s_x)$.

If $\gamma$ is a transformation of $\cX$, then its graph
$\Graph(\gamma)$ 
intersects the fiber $\cS_x$ at a point
$(\gamma(\tilde{x}),\tilde{x})$, 
where $\tilde{x}$ is a fixed point of the mapping $s_x\circ\gamma$. 
If~$\gamma$ is close enough to the identity mapping, then the
fixed point $\tilde{x}$ is close to~$x$ and unique. 
We will also denote~$\tilde{x}$ by~$\tilde{x}^\gamma$ to indicate what
mapping~$\gamma$ generates this fixed point.

The correspondence $x\to \tilde{x}$ is a local diffeomorphism. 
The inverse mapping $\tilde{\gamma}:\, \tilde{x}\to x$ 
we call the {\it mid-trans\-formation\/} related to~$\gamma$.

\begin{theorem}
{\rm(i)} The transformation~$\gamma$ is reconstructed via 
reflections and the mid-trans\-formation~$\tilde{\gamma}$ 
by the formula
\begin{equation}
\gamma(z)=s_{\tilde{\gamma}(z)}(z).
\tag{3.1}
\end{equation}

{\rm(ii)} Let $\gamma$ be symplectic, 
and let $\tilde{x}=\tilde{x}^\gamma$ 
be the fixed point of $s_x\circ\gamma$.
Then the $1$-form
$$
-\cH_x(\tilde{x})\equiv \cH_x(\gamma(\tilde{x}))
$$
is closed. In the simply connected case there is 
a function $\Phi^\gamma$ such that 
\begin{equation}
d\Phi^\gamma(x)+\cH_x(\tilde{x})=0.
\tag{3.2}
\end{equation}

{\rm(iii)} Stationary points of the function $\Phi^\gamma$ 
are fixed points of the symplectic transformation~$\gamma$.
The differential of~$\gamma$ and the matrix
of the second derivatives 
of~$\Phi^\gamma$ at the fixed (stationary) point are related to
each other by the formula
\begin{gather*}
d\gamma=\frac{I+\frac12\Psi\cdot D^2 \Phi^\gamma}
{I-\frac12\Psi\cdot D^2 \Phi^\gamma},
\qquad\text{or}\qquad
D^2 \Phi^\gamma=2\omega \cdot \frac{d\gamma-I}{d\gamma+I}
\\
[\text{at the fixed point}].
\end{gather*}
Here $\Psi=\omega^{-1}$ is the Poisson tensor on~$\cX$. 
\end{theorem}

\begin{proof}
Formula (3.1) is just a consequence of the definition of
$\tilde{x}$ and $\tilde{\gamma}$. 
Assertion (ii) is another statement of the fact 
that~$\gamma$ is symplectic  
(the submanifold $\Graph(\gamma)$ is Lagrangian 
in $\cX\times\cX$).   
Assertion (iii) follows from (3.2) and the boundary
conditions (2.2).
\end{proof}

We call $\Phi^\gamma$ the {\it phase function\/} corresponding
to the transformation~$\gamma$. 
This name is due to the fact that $\Phi^\gamma$ is the phase of
the WKB-approximation of the quantum operator related
to~$\gamma$. 

Note that the phase function is uniquely defined if some 
additional Cauchy data are fixed, say, 
zero data at some point $y\in\cX$. 
We denote such a {\it normalized phase function\/}
corresponding to~$\gamma$ by~$\Phi^\gamma_y$.
Thus, by definition, the function $\Phi^\gamma_y$ obeys (3.2)
and $\Phi^\gamma_y(y)=0$.

Hereafter, 
we assume that the form~$\omega$ is exact.
We need this condition to be sure that integrals of the
symplectic form along membranes, which we consider below, do not
depend on the choice of membrane surface. 
Actually, in applications of this phase analysis 
in quantum theory all the integrals are staying in the exponent
and we can replace the exactness of~$\omega$ by the quantization
condition on the cohomology classes.

\begin{theorem}
{\rm(i)} The normalized phase function of the symplectic
transformation~$\gamma$ is given by the formula
\begin{equation}
\Phi^\gamma_y(x)=\int_{\Sigma^\gamma(x,y)}\omega.
\tag{3.3}
\end{equation}
Here $\Sigma^\gamma(x,y)$ is a membrane in~$\cX$ 
whose boundary is composed by four pieces{\rm:}
an arbitrary curve~$c$ connecting $\tilde{x}$ with $\tilde{y}$, 
the Ether geodesic from~$\tilde{y}$ to $\gamma(\tilde{y})$ 
through the mid-point~$y$, 
the curve~$\gamma(c)$
{\rm(}with the opposite orientation{\rm)}
connecting~$\gamma(\tilde{y})$ 
with~$\gamma(\tilde{x})$,
and the Ether geodesic from~$\gamma(\tilde{x})$ 
to~$\tilde{x}$ through the mid-point~$x$.

{\rm(ii)} The following cocyclic properties hold{\rm:}
\begin{align*}
&\Phi^\gamma_x(y)+\Phi^\gamma_y(x)=0,\qquad 
\Phi^\gamma_y(x)+\Phi^\gamma_w(y)+\Phi^\gamma_x(w)=0,
\\
&\Phi^{\gamma^{-1}}_y(x)+\Phi^\gamma_y(x)=0, \qquad  
\forall x,y,w\in\cX.
\end{align*}
\end{theorem}

\begin{proof}
Assertion (i) can be proved in the same way as Lemma~7.1\,(ii)
in \cite{K};
see also the proof of Theorem~8.1 below. 
Assertion (ii) is a direct consequence of (3.3) and
the Stokes theorem. 
\end{proof}

In the conclusion of this section, 
let us make some remarks
about the case of symmetric symplectic spaces~\cite{BCG}. 

First, note that in the construction of Theorem~3.1  
the main role is played by the Lagrangian fibration 
generated by graphs of reflections. 
The Lagrangiancy is equivalent to symplecticity 
of reflections. In our approach the symplecticity 
is an automatical consequence of the definition of reflections
$s=s_x$ by means of the Hamilton type dynamic equation
\begin{equation}
\frac{\pa}{\pa x}s=D\cH_x(s)\Psi(s),\qquad 
s\bigg|_{x=z}=z.
\tag{3.4} 
\end{equation}
The Ether Hamiltonian $\cH$ is the object which generates 
this dynamics.

So, the symplecticity of~$s_x$ is the priority.
At the same time, we lose (in general) 
the usual relationship of reflections 
with the exponential mappings. 
Our reflections are not geodesic reflections:
$$
s_x\ne s^0_x,\qquad \text{where}\quad
s^0_x(z)\od \exp_x\big(-\exp^{-1}_x(z)\big).
$$
Instead of this, we have
$$
s_x(z)=\Exp_x\big(-\Exp^{-1}_x(z)\big),
$$
where $\Exp_x$ is the Ether exponential mapping.
That is why we use in Theorem~3.2 and everywhere below 
the Ether geodesics, but not the usual $\Gamma$-geodesics. 

In the classical theory of symmetric spaces, 
initiated by E.~Cartan~\cite{EC},
the geodesics and the geodesic reflections $s^0_x$  
play an exclusive role. 
In the symplectic situation~\cite{BCG},
{\it the geodesic reflections 
in order to be symplectic mappings
must satisfy the Loos condition}~\cite{L}
$s^0_x s^0_y s^0_x=s^0_{s^0_x(y)}$ 
or the Cartan condition $\nabla R=0$.
One cannot avoid these conditions 
if only geodesic reflections are considered. 

Under these conditions, i.e., 
in the framework of symmetric symplectic spaces, 
the correspondence between symplectic transformations 
and phase functions was studied in~\cite{Oz}.
The construction in~\cite{Oz} is based on 
the use of exponential mappings $\exp_x$, 
which are the usual ``kinematic'' tools 
in this framework. 
But one can clearly see that, 
besides the local character, the presence 
of exponential mappings in all formulas
actually implies loosing some important geometric structures
like, for instance, the dynamic equation~(3.4) 
or the zero curvature equation~(2.1).
The information carried by 
the Ether Hamiltonian allows one to avoid these kinematic
difficulties and makes the dynamic view~(3.4) to be the basic
point. 
It seems that Eq.~(3.4), 
even for $s=s^0$ in the symmetric case,
was not used in the literature.

\section{Dynamic phase functions}

The most important examples of symplectic transformations are 
translations along trajectories of Hamiltonian systems. 
We denote such a translation by $\gamma^t_H$, 
where $H$ is the Hamilton function and~$t$ is the time variable.
If~$t$ is close to zero, then the corresponding quantum flow in
the WKB-approximation is described by a
{\it dynamic phase function}~$\Phi^t$,
namely,
$$
\exp\bigg\{-\frac{it}{\hb}\wh{H}\bigg\}=\wh{G}^t,
\qquad
G_t=\exp\bigg\{\frac{i}{\hb}\Phi^t\bigg\}\varphi^t+O(\hb)
$$
(in a semiclassically-simple domain); see details in~\cite{K}.

The following formula for $\Phi^t$ was derived in~\cite{K} 
via the membrane area:
\begin{equation}
\Phi^t(x)=\int_{\Sigma^t(x)}\omega-tH(\tilde{x}).
\tag{4.1}
\end{equation}
Here $\Sigma^t(x)$ is a {\it dynamic segment\/} bounded by the
Hamiltonian trajectory (whose time-length is~$t$) and by the Ether
geodesics connecting the ends of the trajectory and passing
through the mid-point~$x$.
The value of the Hamilton function~$H$ in (4.1) is taken on the
trajectory-side of~$\Sigma^t(x)$.

Formula (4.1) is a generalization of the 
mid-point formulas found by Berry and Marinov in Euclidean phase
spaces~\cite{B,M}, see also~\cite{Oz} for symmetric spaces.  
A certain modification using membranes with ``wings'' 
was suggested in~\cite{KO} for magnetic phase spaces.

\begin{theorem}
The function {\rm(4.1)} is the phase function corresponding 
to the Hamiltonian translation~$\gamma^t_H$. 
The normalized phase function of $\gamma^t_H$ 
given by {\rm(3.3)} is related to {\rm(4.1)} 
via the identity{\rm:} 
\begin{equation}
\Phi^{\gamma^t}_y(x)=\Phi^t(x)-\Phi^t(y).
\tag{4.2}
\end{equation}
\end{theorem}

\begin{proof}
The first statement was proved in \cite{K}. 
Formula (4.2) is a consequence of the Stokes theorem 
and the following version of the Poincare--Cartan 
``integral invariant'' formula
\begin{equation}
t\big(H(z)-H(w)\big)=\int_{\Sigma^t_{z,w}}\omega.
\tag{4.3}
\end{equation}
Here the boundary of $\Sigma^t_{z,w}$ 
consists of two pieces of the Hamiltonian trajectories 
passing through~$z$ and~$w$, 
of a path~$c$ connecting~$w$ with~$z$, 
and of the path $\gamma^t_H(c)$.
\end{proof}

\section{Membrane representation}

The name ``phase function'' which we use for solutions of 
Eq.~(3.2) is not the only natural candidate. At the same time one
can use for $\Phi^\gamma$ the name: 
{\it generating function of~$\gamma$}.  
This is due to the following statement.

\begin{theorem}
Let a real smooth function~$\Phi$ be close enough to a constant,
so that the equation
\begin{equation}
d\Phi(x)+\cH_x(z)=0,
\tag{5.1}
\end{equation}
has a smooth solution $x=\tilde{\gamma}(z)$ close enough 
to the identical $x=z$.
Then the transformation~$\gamma$ defined by {\rm(3.1)} 
is symplectic, and~$\tilde\gamma$ is its mid-trans\-formation.
\end{theorem}

Thus the function $\Phi$ generates a symplectic transformation
of~$\cX$. 

\begin{corollary}
Any function $\Phi$ satisfying the conditions of
Theorem~{\rm5.1} can be represented via symplectic area as
follows{\rm:} 
\begin{equation}
\Phi(x)=\int_{\Sigma^\gamma(x,y)}\omega+\Phi(y),\qquad 
\forall x\in\cX.
\tag{5.2}
\end{equation}
Here $\gamma$ is the symplectic transformation generated
by~$\Phi$ as in Theorem~{\rm5.1}, 
and $\Sigma^\gamma(x,y)$ is the membrane in~$\cX$ 
defined in Theorem~{\rm3.2,\,(i)}.
The point $y\in\cX$ is fixed.
\end{corollary}

We call formula (5.2) the {\it membrane representation\/} of
the function~$\Phi$.

\section{Geometric phase product}

The product of quantum operators corresponds to a noncommutative
product of functions over~$\cX$. The integral kernel of the
product operation in semiclassical approximation 
can be described by the phase function
\begin{equation}
\Phi_{y,z}(x)=\int_{\Delta(x,y,z)}\omega.
\tag{6.1}
\end{equation}
Here the membrane $\Delta(x,y,z)$ in $\cX$ is composed by three
Ether geodesics passing through mid-points $z$, $y$, and~$x$
close enough to each other 
(see~\cite{K,W}).

On the level of phase functions the noncommutative product 
is given by 
\begin{equation}
(\Phi''\circ \Phi')(x)\od
\big[\Phi''(x'')+\Phi'(x')
+ \Phi_{x'',x'}(x)\big]_{\substack{x'=X'(x) \\ x''=X''(x)} }\,\,,
\tag{6.2}
\end{equation}
where $X'(x)$ and $X''(x)$ are stationary points of the  
right-hand side of (6.2) with respect to~$x'$ and~$x''$. 

In the case of Euclidean space, 
the phase product~(6.2) was considered in many papers 
dealing with
the Fourier integral operators and the Maslov canonical operator,
for instance, in~\cite{H,GSt,Tr,OMY},
and in the case of symmetric symplectic manifolds,
the detailed study was done in~\cite{Oz}.   

We call (6.2) a {\it phase product\/} over~$\cX$.

\begin{theorem}
{\rm(i)} The phase product is associative and has the unity
element $\Phi=0$.

{\rm(ii)}  The phase product {\rm(6.2)}
of phase functions of symplectic transformations 
{\rm(}close enough to the identity{\rm)} 
is the phase function of a composition of these transformations. 

{\rm(iii)} The family of functions {\rm(4.1)}
forms a one-parameter local group with respect to the
phase product~{\rm(6.2)}. 

{\rm(iv)} The triangle area {\rm(6.1)} 
is the phase {\rm(}or generating{\rm)} function 
corresponding to the Ether translation 
$g_{y,z}=s_y\circ s_z$. 
For this transformation, 
the membrane formula {\rm(3.3)} is equivalent to {\rm(6.1):}
$$
\Phi^{g_{y,z}}_y=\Phi_{y,z}.
$$
\end{theorem}

Now let us take two functions, consider the corresponding
symplectic transformations, 
and make up the phase product of their phase functions.

\begin{theorem}
Let functions $\Phi'$ and $\Phi''$ satisfy the conditions of
Theo-\break rem~{\rm5.1}. 
Let~$\gamma'$ and~$\gamma''$ be the corresponding symplectic
transformations, and let $\gamma=\gamma''\circ \gamma'$
be their composition. 
For any $x\in\cX$ consider the fixed point~$\tilde{x}$ 
of the mapping $s_x\circ\gamma$.
Let~$x'$ be the mid-point of the Ether geodesics
between~$\tilde{x}$ and~$\gamma'(\tilde{x})$,
and let $x''$ be the mid-point of the Ether geodesics 
between~$\gamma'(\tilde{x})$ and~$\gamma(\tilde{x})$.  
Then the phase product of functions~$\Phi'$ and~$\Phi''$ 
is given by  
\begin{equation}
(\Phi''\circ \Phi')(x)
=\Phi''(x'')+\Phi'(x')+ \int_{\Delta(x'',x',x)}\omega.
\tag{6.3}
\end{equation}
The differential of this phase product 
is given by the Ether Hamiltonian{\rm:} 
\begin{equation}
d(\Phi''\circ\Phi')(x)=\cH_x(\gamma(\tilde{x})).
\tag{6.4}
\end{equation}
\end{theorem}

\begin{proof}
By (3.2), we know that
$d\Phi''(x'')=-\cH_{x''}(\widetilde{x''})$,  
where $\widetilde{x''}$ is the fixed point 
of $s_{x''}\circ\gamma''$. 
From~\cite{K}, formula (7.9), we have
$\pa_{x''}\Phi_{x'',x'}(x)=\cH_{x''}(a)$, 
where $a$ is the vertex of $\Delta(x'',x',x)$ 
common for sides with mid-points $x'$ and $x''$.
So, the stationary phase condition in (6.2) implies 
$a=\widetilde{x''}$.
In the same way one can prove that 
$\tilde{x}=\widetilde{x'}$ is the vertex of 
$\Delta(x'',x',x)$ common for sides with mid-points~$x$
and~$x'$. Thus the triple $(X'',X',x)$ in (6.2) coincides 
with the triple $(x'',x',x)$ in Theorem~6.2, and then  
(6.3) follows from (6.2) and~(6.1).
The differential is given by 
$d(\Phi''\circ\Phi')(x)=d\Phi_{x'',x'}(x)=\cH_x(\gamma(\tilde{x}))$,
since $\gamma(\tilde{x})$ is the vertex of $\Delta(x'',x',x)$
common for sides with mid-points~$x$ and~$x''$.
\end{proof}

\begin{corollary}
Let a function $\Phi^0$ be close enough to constant
{\rm(}as in Theorem~{\rm5.1)}, 
and $\Phi^t$ be the dynamic phase function {\rm(4.1)}.
Then the phase product 
$\Phi(x,t)\od (\Phi^t\circ \Phi^0)(x)$ 
is given by 
\begin{equation}
\Phi(x,t)=\Phi^0(x')+\Phi^t(x'')+\int_{\Delta(x'',x',x)}\omega.
\tag{6.5}
\end{equation}
Here $x'$ is defined by $d\Phi^0(x')=\cH_{x'}(a)$, 
where~$a$ is the fixed point 
of the map $s_{x'}\circ s_{x}\circ \gamma^t_H$, 
and~$x''$ is the mid-point of the Ether geodesic from 
$a$ to $\gamma^t_H (a)$.
\end{corollary}

One can call formulas like (6.3) and (6.5)
the {\it geometric representation of the phase product}.

In particular, consider the family 
of dynamic phase functions $\Phi^t$~(4.1). 
On the membrane level we have
\begin{equation}
\Sigma^{\tau+t}(x)=\Sigma^\tau(X^\tau)\cup\Sigma^t(X^t)
\cup\Delta(X^\tau,X^t,x).
\tag{6.6}
\end{equation}
Here $X^t$ is the mid-point of the Ether geodesics 
between~$\tilde{x}$ and~$\gamma^t(\tilde{x})$, 
$X^\tau$ is the mid-point between~$\gamma^t(\tilde{x})$ 
and~$\gamma^{\tau+t}(\tilde{x})$,
and we denote by $\tilde{x}=\tilde{x}^{\gamma^{\tau+t}}$
the fixed point of $s_x\circ\gamma^{\tau+t}$.

Thus by formula~(4.1)
we have
\begin{align*}
\Phi^{\tau+t}(x)
&=\int_{\Sigma^{\tau+t}(x)}\omega-(\tau+t)H(\tilde{x})
\\[3\jot]
&\!\!\!\!\!\!\overset{\text{see (6.6)}}{=}
\bigg(\int_{\Sigma^\tau(X^\tau)}\omega-\tau H\bigg)
+\bigg(\int_{\Sigma^t(X^t)}\omega-tH\bigg)
+\int_{\Delta(X^\tau,X^t,x)}\omega
\\[3\jot]
&\!\!\!\!\!\!\overset{\text{see (6.3)}}{=}
(\Phi^\tau\circ\Phi^t)(x).
\end{align*}

{\it Formula} (6.6) {\it and the last calculation geometrically 
represent the statement of Theorem~{\rm6.1,\,(iii):} 
the group property of the family of phase functions 
$\Phi^t$ with respect to the phase product on general symplectic
manifolds}. 

In Euclidean spaces this calculation was first
demonstrated in~\cite{M}. 

Formula (6.5) can also be easily interpreted in this way.
Actually, there is a natural extension of this representation 
for phase functions of general symplectic transformations.

Let $(x'',x',x)$ be a triple of points related to symplectic 
transformations $\gamma'$ and $\gamma''$ as in Theorem~6.2, and
let $(y'',y',y)$ be another such triple. Then, on the membrane
level, we have 
\begin{equation}
\Sigma^{\gamma''\circ \gamma'}(x,y)\cup \Delta(y'',y',y)
=\Sigma^{\gamma''}(x'',y'')\cup\Sigma^{\gamma'}(x',y')
\cup \Delta(x'',x',x).
\tag{6.7}
\end{equation}

By applying the Stokes theorem and formula (6.3), one obtains:
\begin{equation}
\Phi^{\gamma''}_{y''}\circ\Phi^{\gamma'}_{y'}
=\Phi^{\gamma''\circ \gamma'}_y +\Phi_{y'',y'}(y).
\tag{6.8}
\end{equation}
Here, on the left, we have the phase product of normalized
generating functions of two symplectic transformations, 
and, on the right, we have a generating function of the
composition of these transformations.

The normalized generating functions are given by the membrane
area via~(3.3). Thus {\it formulas} (6.7), (6.8) 
{\it represent the statement of Theo-\break rem~{\rm6.1,\,(ii)}
geometrically via symplectic areas of membranes on general
symplectic manifolds}.

\section{Groupoid interpretation of the phase\\ product}

The boundary condition (2.2) guarantees that 
the matrix $D\cH_x(z)$ is not degenerate  
at the diagonal $\{x=z\}\subset \cX\times\cX$.
Denote by $\cX^\#\subset \cX\times\cX$ 
a connected reflective neighborhood of the diagonal 
where $\det D\cH_x(z)\ne0$.
Also denote by $\cE$ the corresponding neighborhood 
of the zero section in $T^*\cX$:
$$
\cE=\{(x,p)\in T^*\cX\mid p=\cH_x(z),\,(x,z)\in\cX^\# \}.
$$
Then there is a fibration of~$\cE$ over~$\cX$,
\begin{equation}
\ell:\, \cE\to\cX,\qquad 
\ell(x,p)\od z\qquad\text{if}\quad p=\cH_x(z),
\tag{7.1}
\end{equation}
and a dual fibration 
\begin{equation}
r:\, \cE\to\cX,\qquad 
r(x,p)\od \ell(x,-p).
\tag{7.1\,a}
\end{equation}

In view of (2.3), we have
\begin{equation}
\ell(x,p)=s_x(r(x,p)).
\tag{7.2}
\end{equation}
The zero curvature equation (2.1) implies that~$\ell$ 
is a Poissonian mapping and~$r$ is an anti-Poissonian mapping
commuting with~$\ell$, that is, 
\begin{equation}
\{\ell^j,\ell^k\}=\Psi^{jk}(\ell),\qquad 
\{r^j,r^k\}=\Psi^{kj}(\ell),\qquad 
\{\ell^j,r^k\}=0,
\tag{7.3}
\end{equation}
where the brackets $\{\cdot,\cdot\}$ correspond to the standard
symplectic form $dp\wedge dx$ on~$\cE$.
System (7.3) is known as the Lie--Engel system.

Thus we conclude that in the sense of \cite{KM1,KM2}
the space~$\cE$ is a {\it phase space\/} over~$\cX$
equipped with the Poisson bifibration (7.1), (7.1\,a), and
(7.3). 
Actually, such bifibrations 
(in the case of Poisson brackets of constant rank)
were first considered by S.~Lie; 
see references for the general Poisson case in~\cite{KM1}
and for more details in \cite{KM2}--\cite{KMbook}.

The boundary conditions (2.2) imply 
\begin{gather}
\ell(x,p)\bigg|_{p=0}=x,\qquad 
\frac{\pa \ell(x,p)}{\pa p}\bigg|_{p=0}=\frac12\Psi(x),
\tag{7.4}\\
\frac{\pa^2 \ell(x,p)}{\pa p\pa p}\bigg|_{p=0}
=\frac14\Psi(x)\Gamma(x)\Psi(x).
\nonumber
\end{gather}
Eqs.~(7.2) and~(7.4) relate the phase space structure on~$\cE$
with the reflective structure on~$\cX$, in particular, 
with the symplectic connection~$\Gamma$ on~$\cX$
(see~\cite{K}).

Each function~$H$ on~$\cX$ is lifted up to 
the function $\ell^*H$ on $\cE\subset T^*\cX$.
The last one can be considered as a Hamilton function 
for the Hamilton--Jacobi equation over~$\cX$ 
(see in~\cite{KM1} for the general Poisson case):
\begin{equation}
\frac{\pa \Phi}{\pa t}
+H\bigg(\ell\bigg(x,\frac{\pa\Phi}{\pa x}\bigg)\bigg)=0.
\tag{7.5}
\end{equation}

\begin{theorem}
The dynamic phase function $\Phi^t$ {\rm(4.1)}
is the solution of the Hamilton--Jacobi equation~{\rm(7.5)}
with zero Cauchy data. The phase product function~$\Phi$
{\rm(6.5)} is the solution of~{\rm(7.5)}
with the Cauchy data $\Phi\big|_{t=0}=\Phi^0$.
\end{theorem}

\begin{proof}
The first statement was proved in~\cite{K}.
From formula (6.3) we see that  
\begin{equation}
\frac{\pa}{\pa t}\Phi(x,t)
=\frac{\pa \Phi^t}{\pa t}(x')
=-H\big(\ell(x',d\Phi^t(x'))\big).
\tag{7.6}
\end{equation}

As in the proof of Theorem~6.2, we know that 
$d\Phi^t(x')=\cH_{x'}(\gamma^t(a))$
and thus, 
$\ell(x',d\Phi^t(x'))=\gamma^t_H(a)=\gamma^t_H(\gamma^0(\tilde{x}))$,
where~$\gamma^0$ is the symplectic transformation generated 
by the function~$\Phi^0$. 

From (6.4) we have
$$
d_x\Phi(x,t)=d(\Phi^t\circ \Phi^0)(x)
=\cH_x\big(\gamma^t_H(\gamma^0(\tilde{x}))\big)
$$
and thus 
$\ell(x,d_x\Phi(x,t))=\gamma^t_H(\gamma^0(\tilde{x}))$.

So we see that 
$\ell(x',d\Phi^t(x'))=\ell(x,d_x\Phi(x,t))$. 
Then Eq.~(7.5) follows from~(7.6).
\end{proof}

Let us now consider the phase product (6.3) from the groupoid
point of view.

The phase space $\cE\subset T^*\cX$ can be represented 
in a neighborhood of the diagonal in $\cX\times\cX$ 
by means of the Poisson bifibration
\begin{gather*}
\ell\times r:\, \cE\to \cX\times\cX^{(-)},\qquad
(\ell\times r)(x,p)\od \big(\ell(x,p),r(x,p)\big),
\\
(\ell\times r)^* dp\wedge dx=\omega(\ell)-\omega(r).
\end{gather*}

In the direct product $\cX\times\cX$,  
there is a natural groupoid multiplication
\begin{equation}
(x,y)\otimes(y,z)=(x,z).
\tag{7.7}
\end{equation}

Transporting this multiplication back to~$\cE$
by means of $(l\times r)^{-1}$, 
we obtain the groupoid structure
\begin{equation}
(x,p)\circledcirc(y,\xi)\od (z,\eta)
\tag{7.8}
\end{equation}
just by solving the system of equations
\begin{gather*}
\ell(z,\eta)=\ell(x,p),\qquad r(z,\eta)=r(y,\xi),
\\
r(x,p)=\ell(y,\xi).
\end{gather*}
Here $p\in T^*_x \cX$, $\xi\in T^*_y\cX$, $\eta\in T^*_z\cX$.
In the symplectic case that we consider, 
the multiplication (7.8) coincides with that 
described in~\cite{K1} 
for the case of general Poisson manifolds.

With respect to the multiplication (7.8),
the mappings~$\ell$ and~$r$ are just 
the left and right (or the source and target) groupoid mappings:
\begin{equation}
\ell(m)=m\circledcirc m^{-1},\qquad
r(m)=m^{-1}\circledcirc m,\qquad
m\in\cE.
\tag{7.9}
\end{equation}
There is a consistency condition between the groupoid
multiplication (7.8) 
and the symplectic structure in~$\cE$;
namely, the mapping~${\circledcirc}$
preserves the symplectic structure 
or the graph of $\circledcirc$ in a Lagrangian submanifold 
in $\cE\times\cE\times\cE^{(-)}$,
see details in~\cite{K1,KMbook,W2}.

For any two subsets $\Lambda',\Lambda''\subset\cE$,
one can define their groupoid product: 
\begin{equation}
\Lambda''\circledcirc\Lambda'
\od 
\{m''\circledcirc m'\mid m''\in\Lambda'',\,m'\in\Lambda'\}.
\tag{7.10}
\end{equation}

\begin{lemma}
{\rm(i)}
Let $\Lambda'$ and $\Lambda''$ be submanifolds, 
the left mapping~$\ell$ {\rm(7.9)} restricted to~$\Lambda'$
be a diffeomorphism, and the right mapping~$r$ {\rm(7.9)} 
restricted to $\Lambda''$ be a diffeomorphism.
Then the subset {\rm(7.10)} is a submanifolds.

{\rm(ii)}
If $\Lambda''$ and $\Lambda'$ are Lagrangian, 
then $\Lambda''\circledcirc\Lambda'$ is Lagrangian
at every point where it is a submanifold.
\end{lemma}

Now, to each function $\Phi$ on~$\cX$ one can assign 
a Lagrangian submanifold in~$T^*\cX$:
\begin{equation}
\Lambda^\Phi=\{(x,p)\mid p=d\Phi(x)\}.
\tag{7.11}
\end{equation}
If the function~$\Phi$ is close to constant, 
then the submanifold $\Lambda^\Phi$ belongs to the groupoid 
$\cE\subset T^*\cX$.
We call $\Phi$ the {\it generating function\/}
of~$\Lambda^\Phi$. 

Let us take two functions~$\Phi'$ and $\Phi''$ 
that satisfy conditions of Theorem~5.1 
and consider their phase product $\Phi'\circ\Phi''$  
as in Theorem~6.2.
The following statement is just a reformulation
of the construction described in Theorem~6.2.

\begin{theorem}
The phase product {\rm(6.2)} of functions over~$\cX$
corresponds to the groupoid product {\rm(7.10)}
of Lagrangian submanifolds in $\cE\subset T^*\cX$, 
that is,
\begin{gather}
\Lambda^{\Phi''} \circledcirc \Lambda^{\Phi'}
= \Lambda^{\Phi''\circ \Phi'},
\nonumber\\
\Lambda^{\Phi} \circledcirc \Lambda^{-\Phi}
\subset \Lambda^{0}\equiv\cX,
\tag{7.12}\\
\Lambda^{\Phi} \circledcirc \Lambda^{0}
=\Lambda^{0} \circledcirc \Lambda^{\Phi}
=\Lambda^{\Phi}.
\nonumber
\end{gather}
Here the submanifold $\Lambda^{0}=\cX$
{\rm(}the zero section in $T^*\cX${\rm)}
is assigned to the zero function $\Phi=0$,
and so the base manifold~$\cX$ plays the role 
of the unity element for the product~{\rm(7.10)}.
\end{theorem}

\section{Chord submanifolds}

Actually, Theorem~7.3 is a consequence of the general
observation~\cite{K1} 
that on the quantum level 
the groupoid product~$\circledcirc$ corresponds 
to a noncommutative algebra of operators (quantum observables). 
This algebra admits very tricky constructions.
Thus we can expect that the groupoid product might be used 
to create some interesting extensions of objects 
of symplectic geometry.

The simplest construction of such a type is inspired by the 
known quantum Weyl--Wigner isomorphism between integral kernels
and symbols of operators.
In the symplectic geometry this isomorphism is represented by
the left-right groupoid mapping 
$\ell\times r:\, \cE\to \cX\times\cX^{(-)}$.

To each Lagrangian submanifold $M\subset\cX\times \cX^{(-)}$ 
(on which the symplectic form $\omega\circleddash\omega$ is
annulated), 
one can assign a Lagrangian submanifold in~$\cE$:
$$
\Lambda_M\od\{m\in\cE\mid 
(m\circledcirc m^{-1},m^{-1}\circledcirc m)\in M\}
=(\ell\times r)^{-1}(M).
$$

In particular, if $M=\Graph(\gamma)$, 
where $\gamma$ is a symplectic transformation of~$\cX$, then 
$$
\Lambda_{\Graph(\gamma)}
=\{m\in\cE\mid m\circledcirc m^{-1}=\gamma(m^{-1}\circledcirc m)\}.
$$
Assuming that $\gamma$ is close enough to the identical mapping,
we can represent this Lagrangian submanifold in the form~(7.10),
that is, 
$$
\Lambda_{\Graph(\gamma)}=\Lambda^{\Phi^{\gamma}},
$$
where $\Phi^{\gamma}$ is the generating function corresponding
to~$\gamma$ via Theorem~5.1. 

Another important case is $M=\lambda\times\lambda$, 
where~$\lambda$ is a Lagrangian submanifold in~$\cX$.
In this case 
\begin{equation}
\Lambda_{\lambda\times\lambda}
=\{m\in\cE\mid m\circledcirc m^{-1}\in\lambda,\,
m^{-1}\circledcirc m\in\lambda\}.
\tag{8.1}
\end{equation}

For each $m=(x,p)\in\Lambda_{\lambda\times\lambda}$ 
we have two points $b\od r(x,p)$ and $a\od\ell(x,p)$ 
belonging to~$\lambda$.
One can identify~$m$ with the Ether geodesic connecting~$b$
with~$a$ and passing through the mid-point~$x$.
Such a geodesic can be called a {\it chord\/}
of the submanifold~$\lambda$.
So $\Lambda_{\lambda\times\lambda}$ can be considered 
as a set of all chords. 
We call $\Lambda_{\lambda\times\lambda}$
a {\it chord submanifold}.

Note that here we have in mind the oriented chords, 
but one can consider the nonoriented chords as well.

\begin{theorem}
Let $\wt{\cX}_\lambda\subset \cX$ be a connected domain 
such that for any $x\in\wt{\cX}_\lambda$ 
there is a unique nonoriented chord of~$\lambda$ passing 
through the mid-point~$x$.

Then over the domain $\wt{\cX}_\lambda$ 
the chord Lagrangian submanifold
$\Lambda_{\lambda\times\lambda}$ {\rm(8.1)} 
can be represented in the form analogous to {\rm(7.10):}
$$
\Lambda_{\lambda\times\lambda}=\{(x,p)\mid p=\pm d\Phi_\lambda(x)\},
$$
where the function $\Phi_\lambda$ is given by the integral
\begin{equation}
\Phi_\lambda(x)=\int_{\Sigma_\lambda(x)}\omega.
\tag{8.2}
\end{equation}
Here $\Sigma_\lambda(x)$ is a membrane in~$\cX$ composed of the
chord and the end-points of this chord.
The differential of $\Phi_\lambda$ is given by 
\begin{equation}
d\Phi_\lambda(x)=\cH_x(a),
\tag{8.3}
\end{equation}
where $a$ is the end-point of the chord.
\end{theorem}

\begin{proof}
We need to prove that the end-points $a,b\in\lambda$ 
of the Ether geodesics passing through the mid-point 
$x\in\wt{\cX}_\lambda$ are related to the function~(8.2) 
by means of the equations
\begin{equation}
a=\ell(x,d\Phi_\lambda(x)),\qquad b=r(x,d\Phi_\lambda(x)).
\tag{8.4}
\end{equation}

Let us fix a certain point $x_0\in\wt{\cX}_\lambda$ close enough
to~$x$ and denote by $a_0,b_0\in\lambda$ 
the end-points of the corresponding Ether geodesic through 
the mid-point~$x_0$.
The function $\Phi_\lambda(x)$ (8.2) can be represented as 
\begin{equation}
\Phi_\lambda(x)=\Phi_\lambda(x_0)+\int_{\Sigma_\lambda(x,x_0)}\omega.
\tag{8.5}
\end{equation}
Here $\Sigma_\lambda(x,x_0)$ is a slice between two fixed
Ether geodesics passing through the mid-points~$x$ and~$x_0$. 
The part of the boundary of $\Sigma_\lambda(x,x_0)$ 
belonging to~$\lambda$ consists of two paths 
$A=\{A(t)\mid t\in [0,1]\}$ and $B=\{B(t)\mid t\in [0,1]\}$, 
so that  
$A(0)=a_0$, $A(1)=a$, and $B(0)=b_0$, $B(1)=b$.
The points $A(t)$ and $B(t)$ are just the left and right ends 
of the Ether geodesics fibrating the slice
$\Sigma_\lambda(x,x_0)$.  
Let $X=\{X(t)\mid t\in [0,1]\}$ be the corresponding curve of
mid-points of those Ether geodesics, so that 
$X(0)=x_0$, $X(1)=x$.

Then the curve $X$ separates the slice $\Sigma_\lambda(x,x_0)$ 
in two segments (the left and the right):
\begin{equation}
\Sigma_\lambda(x,x_0)
=\Sigma^{\rm left}_\lambda\cup \Sigma^{\rm right}_\lambda.
\tag{8.6}
\end{equation}
The left segment $\Sigma^{\rm left}_\lambda$ is bounded 
by the left ends $A(t)$ of the Ether geodesics 
and the right segment $\Sigma^{\rm right}_\lambda$ is bounded 
by the right ends $B(t)$. 

Note that both segments are images of one and the same
``vertical'' membrane $\sigma\subset T^*\wt{\cX}_\lambda$, 
that is  
\begin{equation}
\Sigma^{\rm left}_\lambda=\ell(\sigma),\qquad
\Sigma^{\rm right}_\lambda=r(\sigma^{(-)}),
\tag{8.7}
\end{equation}
where the sign minus marks the inversion of the orientation. 
The boundary of the membrane~$\sigma$ is composed of four
pieces. 
The part of the boundary belonging to~$\wt{\cX}_\lambda$ is the
curve~$X$. 
The vertical part of the boundary of~$\sigma$
consists of two paths in the vertical fibers 
$T^*_{x_0}\wt{\cX}_\lambda$ and $T^*_{x}\wt{\cX}_\lambda$.
Each vertical path is projected by the mappings~$\ell$ and~$r$ 
onto the left and right parts of the Ether geodesics 
through~$x_0$ and~$x$. 
The last piece of the boundary of~$\sigma$ is a curve 
$m=\{(X(t),P(t))\mid t\in[0,1]\}
\subset\Lambda_{\lambda\times\lambda}$. 

The points of this curve are mapped by the left and right
mappings to the left and right pieces of the boundary 
of $\Sigma_\lambda(x,x_0)$, 
that is,
\begin{equation}
\ell(X(t),P(t))=A(t),\qquad r(X(t),P(t))=B(t).
\tag{8.8}
\end{equation}

In view of (8.6) and (8.7), we have
\begin{align}
\int_{\Sigma_\lambda(x,x_0)}\omega
&=\int_{\Sigma^{\rm left}_\lambda}\omega
+\int_{\Sigma^{\rm right}_\lambda}\omega
=\int_{\ell(\sigma)}\omega-\int_{r(\sigma)}\omega
=\int_{\sigma}(\ell\times r)^*(\omega\circleddash\omega)
\nonumber\\
&=\int_{\sigma} dp\wedge dx
=\int_{m} pdx 
=\int^1_0 P(t)\,dX(t).
\tag{8.9}
\end{align}

Using formula (8.5), we conclude that 
\begin{equation}
d\Phi_\lambda(X(t))=P(t)=\cH_{X(t)}(A(t)).
\tag{8.10}
\end{equation}
Now taking into account (8.8) and setting $t=1$,
we obtain (8.4) and (8.3).
\end{proof}

Note that one may use a description of the Lagrangian
submanifold $\lambda\subset\cX$ as the joint energy level of
Hamiltonians in involution (on~$\lambda$), namely, 
$$
\lambda=\{x\in\cX\mid H_1(x)=E_1,\dots,H_n(x)=E_n\},
$$
where $E_j$ are constants.
Then the chord function $\Phi_\lambda$ (8.2) can be considered
as a solution of the system of Hamilton--Jacobi equations
\begin{equation}
H_j\big(\ell(x,d\Phi_\lambda(x))\big)
=H_j\big(r(x,d\Phi_\lambda(x))\big)=E_j,
\qquad j=1,2,\dots,n.
\tag{8.11}
\end{equation}

Indeed, this fact follows from (8.4) if one notes that 
$H_j(a)=H_j(b)=E_j$ since $a,b\in\lambda$.

Of course, the function $-\Phi_\lambda$ also satisfies system
(8.11). 

Actually, the sequence of statements (8.5), (8.9), (8.10) 
can be revised to obtain the following result.

\begin{corollary}
In the domain $\wt{\cX}_\lambda$ the chord function
$\Phi_\lambda$ {\rm(8.2)} is a unique {\rm(}up to the sign{\rm)}
solution of the system of Hamilton--Jacobi equations  
{\rm(8.11)} obeying the boundary condition 
$\Phi_\lambda\big|_{\lambda}=0$.
\end{corollary}

In the Euclidean $2$-dimensional case $\cX=\bR^2$, 
the membrane formula (8.2) for the WKB-phase of the Wigner
function was established by M.~Berry~\cite{B} 
in the framework of the semiclassical approximation theory. 
In this case, the Ether geodesics in the definition of the
membrane $\Sigma_{\lambda}(x)$ in (8.2) is just a straight chord
connecting a pair of points of the curve $\lambda\subset\bR^2$.

\begin{remark}
Formula (8.2) is easily generalized 
to the case of fibrated co\-isotropic submanifolds.
Namely, let $\lambda\subset\cX$ be coisotropic and its isotropic
foliation be a fibration (see~\cite{GSt1}).
Denote by $\lambda^\#$ the Whitney sum of two copies of~$\lambda$ 
with respect to this isotropic fibration.
Then~$\lambda^\#$ is a Lagrangian submanifold in $\cX\times\cX^{(-)}$
and we can assign to it the Lagrangian submanifold 
\begin{align*}
\Lambda_{\lambda^\#}
&\od\{m\in\cE\mid
\text{$m\circledcirc m^{-1}$ and $m^{-1}\circledcirc m$ belong}\\
&\qquad
\text{to one and the same isotropic fiber in~$\lambda$}\}. 
\end{align*}
Then in a certain domain $\wt{\cX}_\lambda$, one can represent
this submanifold as 
$$
\Lambda_{\lambda^\#}\od\{(x,p)\mid p=\pm d\Phi_\lambda(x)\}.
$$

{\it Formula {\rm(8.2)} works in this fibrated isotropic case as
well, but in the construction of the membrane
$\Sigma_\lambda(x)$ one has to consider only those paths
on~$\lambda$ whose ends belong to one and the same fiber.}
\end{remark}

We conclude this section with an application 
of the groupoid product construction (7.10) 
and the phase product construction (6.3) 
to the case of chord submanifolds. 

\begin{theorem}
{\rm(i)} Let $\lambda\subset\cX$ be a Lagrangian submanifold,
and let $\gamma:\,\cX\to\cX$ be a symplectic mapping. 
Then 
$$
\Lambda_{\Graph(\gamma)}\circledcirc\Lambda_{\lambda\times\lambda}
=\Lambda_{\gamma(\lambda)\times\lambda}.
$$

Let $\wt{\cX}_\gamma\subset \cX$ be a connected domain such that
for any $x\in \wt{\cX}_\gamma$
there is a unique Ether geodesic through the mid-point~$x$
connecting~$\lambda$ with $\gamma(\lambda)$. 
Let $y\in\wt{\cX}_\gamma$ be a point such that 
$\tilde{y}^\gamma\in\lambda$ {\rm(}see Sec.~{\rm3)}.
Then the phase product of the generating function of~$\gamma$
and the chord function {\rm(8.2)} is given by 
\begin{equation}
(\Phi^\gamma_y\circ\Phi_\lambda)(x)
=\int_{\Sigma^\gamma_\lambda(x,y)}\omega.
\tag{8.12}
\end{equation}
Here the boundary of the membrane 
$\Sigma^\gamma_\lambda(x,y)$ consists of two Ether
geodesics through the mid-points~$x$ and~$y$ 
connecting~$\lambda$ with~$\gamma(\lambda)$ 
and of two paths on~$\lambda$ and on~$\gamma(\lambda)$ 
connecting the end-points of those Ether geodesics.

{\rm(ii)} Let $\Phi^t$ be the dynamic phase function {\rm(4.1)} 
corresponding to a Hamilton flow $\gamma^t_H$.
Then over the domain $\wt{\cX}_{\gamma^t_H}\subset\cX$ 
the phase product of $\Phi^t$ with the chord function {\rm(8.2)}
is given by 
\begin{equation}
(\Phi^t\circ\Phi_\lambda)(x)
=\int_{\Sigma^t_\lambda(x)}\omega-Ht.
\tag{8.13}
\end{equation}
Here the boundary of the membrane 
$\Sigma^t_\lambda(x,y)$ is composed by the Ether
geodesic through the mid-point~$x$ 
connecting~$\lambda$ with~$\gamma^t_H(\lambda)$,
by a Hamiltonian trajectory 
{\rm(}whose time length is~$t${\rm)} 
coming from~$\lambda$ to~$\gamma^t_H(\lambda)$,
and by a path on~$\lambda$
connecting the origin  of this trajectory 
with the end-points of those Ether geodesics.
The function {\rm(8.13)} is the solution 
of the Hamilton--Jacobi equation {\rm(7.5)}
with the initial data $\Phi\big|_{t=0}=\Phi_\lambda$.
\end{theorem}

\section{Groupoid extension of Lagrangian\\ submanifolds}

The groupoid multiplication provides an extension of
submanifolds. For any $\Lambda\subset\cE$ we defined 
the extension $\Lambda^\#\subset\cE\times\cE$ as follows 
\begin{equation}
\Lambda^\#\od\{(m'',m')\mid m''\circledcirc m'\in\Lambda\}.
\tag{9.1}
\end{equation}
Thus $\Lambda^\#$ consists of those pairs of multiplicable
elements of the groupoid~$\cE$ whose product belongs
to~$\Lambda$. 

\begin{lemma}
{\rm(i)}
If $\Lambda$ is Lagrangian, 
then $\Lambda^\#$ is Lagrangian 
{\rm(}at all points where it is a submanifold{\rm)}.

{\rm(ii)}
If $\Phi$ is the generating function of $\Lambda\subset\cE$ 
in the sense of {\rm(7.11)},
then the generating function of 
$\Lambda^\#\subset \cE\times\cE$ is given by
\begin{equation}
\Phi^\#(x,y)=[\Phi(z)+\Phi_{y,x}(z)]\bigg|_{z=Z(x,y)},
\tag{9.2}
\end{equation}
where $Z(x,y)$ is the stationary point of the right-hand side 
of {\rm(9.2)} with respect to~$z$.

{\rm(iii)}
Let $\Lambda=\Lambda_M$, where $M$ is a Lagrangian submanifold
in $\cX\times\cX$. 
For each pair $(x,y)\in\cX\times\cX$ close to the diagonal,
let us denote by $(a,b)$ the point of intersection of~$M$ 
with the graph of the Ether translation 
$s_x\circ s_y$.
Then the stationary point $Z(x,y)$ in {\rm(9.2)} coincides with the
mid-point of the Ether geodesic connecting~$a$ and~$b$.
The generating function $\Phi^\#_M$ of $\Lambda^\#_M$ 
obeys the equations
\begin{equation}
\pa_x \Phi^\#_M(x,y)=\cH_x(b),\qquad
\pa_y \Phi^\#_M(x,y)=-\cH_y(a).
\tag{9.3}
\end{equation}

{\rm(iv)}
If $M=\lambda\times\lambda$, where $\lambda$ 
is a Lagrangian submanifold in~$\cX$, then 
\begin{equation}
\Phi^\#_{\lambda\times\lambda}(x,y)
=\int_{\Sigma_\lambda(Z(x,y))}\omega
+\int_{\Delta(Z(x,y),y,x)}\omega,
\tag{9.4}
\end{equation}
where the membrane $\Sigma_\lambda$ is defined 
in Theorem~{\rm8.1}, and~$\Delta$ is the triangle 
from {\rm(6.1)}.

{\rm(v)}
Let $M=\Graph(\gamma^t_H)$, where $\gamma^t_H$ 
is a Hamiltonian flow in~$\cX$.  
Then the generating function $\Phi^\#_{\Graph(\gamma^t_H)}$ 
is given by a formula similar to {\rm(9.4)} 
with the first summand replaced by the dynamic phase function 
$\Phi^t(Z(x,y))$ given by {\rm(4.1)}.

{\rm(vi)}
Let $y=Y(x)$ be a point such that the pair 
$(s_x(y),y)$ belongs to $(\ell\times r)(\Lambda)$. 
Then $y=Y(x)$ is the stationary point of 
$\Phi^\#(x,y)$ with respect to~$y$, and 
\begin{equation}
\Phi^\#(x,y)\bigg|_{y=Y(x)}=\Phi(x).
\tag{9.5}
\end{equation}
\end{lemma}

Now let us note that submanifolds in $\cE\times\cE$ can be
considered as ``operators'' acting in the space of submanifolds
in~$\cE$. 
Namely, if $A\subset \cE\times\cE$ and $L\subset\cE$, then 
\begin{equation}
A(L)\od\{m\in\cE\mid \exists \wt{m}\in \cE:\,
(m,\wt{m})\in A,\, \wt{m}^{-1}\in L\}.
\tag{9.6}
\end{equation}

The role of the unity operator is played by the submanifold
\begin{equation}
\cX^\#\od\{(m,m^{-1})\mid m\in\cE \}.
\tag{9.7}
\end{equation}
Indeed, $\cX^\#(L)=L$ for any $L\subset\cE$.

The composition of two ``operators'' 
$A,B\subset\cE\times\cE$ is defined as
\begin{equation}
A\otimes B\od\{(m'',m')\mid \exists m:\,(m'',m)\in A,\,
(m^{-1},m')\in B\}.
\tag{9.8}
\end{equation}

Obviously, this is an associative product, 
and the element (9.7) is the unity element 
with respect to this product:
$$
A\otimes \cX^\#=\cX^\#\otimes A=A.
$$

There is a natural inverse
$$
A^{-1}\od\{(m,\wt{m})\mid (\wt{m}^{-1},m^{-1})\in A\}.
$$
It is easy to check that 
$$
A\otimes A^{-1}\subset \cX^\#\qquad
(\text{or}\quad A^{-1}\otimes A\subset \cX^\#)
$$
if $A$ is one-to-one projected to the right
(or the left) multiplier in $\cE\times\cE$.

The product (9.8) is consistent with the action (9.6), that is,
$$
A(B(L))=(A\otimes B)(L).
$$
The permutation operation $A\to A'$,
$$
A'\od\{(m,\wt{m})\mid (\wt{m},m)\in A\},
$$
is also consistent with the product (9.8):
$$
(A\otimes B)'=B'\otimes A'.
$$
For each $\Lambda\subset\cE$, 
we denote $\Lambda^\&\od (\Lambda^\#)'$.

\begin{theorem}
The following properties hold{\rm:}
\begin{align*}
\Lambda^\#_1(\Lambda_2)&=\Lambda_1\circledcirc \Lambda_2,
\\
\Lambda^\&_1(\Lambda_2)&=\Lambda_2\circledcirc \Lambda_1,
\\
\Lambda^\#_1\otimes \Lambda^\#_2 &=(\Lambda_1\circledcirc \Lambda_2)^\#,
\\
\Lambda^\&_1\otimes \Lambda^\&_2 &=(\Lambda_2\circledcirc \Lambda_1)^\&,
\\
\Lambda^\#_1\otimes \Lambda^\&_2 &=\Lambda^\&_2\otimes\Lambda^\#_1.
\end{align*}
\end{theorem}

These properties repeat, on the level of Lagrangian submanifolds,
the properties the left and right quantum left mappings
related to the $*$-product operation over~$\cX$
(see \cite{K,KMbook}).

\section{Generalization to the torsion case}

The symplectic groupoid approach described in Secs.~7--9
is mostly based on the Poisson bifibration (7.1)--(7.3).
In definition (7.1), we used the Ether Hamiltonian~$\cH$ 
which obeys the zero curvature condition (2.1)
and the additional boundary and skew-symmetry conditions (2.2),
(2.3). 
These additional conditions actually are not necessary.
One can generalize them, but still keep 
the appropriate properties of the basic dynamic equation (3.4)
for the family of symplectic mappings $\{s_x\}$. 

One can keep the fixed point property (2.4), 
but refuse the involution property (2.5). 
In this case, the family $\{s_x\}$ will still generate 
a symplectic connection on~$\cX$, but in a more complicated way
than via the simplest formula (2.6), and this connection 
will no longer be torsion free.

The goal in this section is to consider this torsion case and to
show at which places the differences 
from the torsion free case do appear.  

The basic zero curvature equation (2.1) 
for $\cH$ is kept unchanged, as well as the first boundary
condition (2.2), that is
$$
\cH\big|_{\rm diag}=0.
$$
Under these general conditions 
we call $\cH$ an {\it internal Hamiltonian\/} over~$\cX$.
The family of symplectic mappings $\{s_x\}$ is defined by (3.4),
so that condition (2.4) still holds. 
We call $s_x$ {\it inversions}.

The connection on~$\cX$ generated by the family $\{s_x\}$ 
is given by
\begin{equation}
\Gamma^l_{jk}(x)=-\frac{\pa^2 s^l_x(z)}{\pa z^m \pa x^r}
\bigg[\frac{\pa s_x(z)}{\pa z}\bigg]^{-1\,m}_{\,\,k}
\bigg[\frac{\pa s_x(z)}{\pa x}\bigg]^{-1\,r}_{\,\,j}
\bigg|_{z=x},
\tag{10.1}
\end{equation}

\begin{lemma}
The connection $\Gamma$ {\rm(10.1)} is symplectic.
\end{lemma}

In general, this connection is not torsion free.

The second and the third boundary conditions in (2.2), 
as well as the skew-symmetry condition (2.3), 
do not hold for the internal Hamiltonian $\cH$, 
in general.
But still the following boundary condition holds:
$$
\nabla_l\nabla_m\cH_k\bigg|_{\rm diag}
+\nabla_s\cH_k\bigg|_{\rm diag} \Psi^{sr} T^j_{rl}\omega_{jm}=0,
$$
where $\nabla$ is taken with respect to the connection~$\Gamma$ 
(10.1), and~$T$ is the torsion tensor of~$\Gamma$.

Now, the left mapping $\ell$ of the bifibration  of the phase
space over~$\cX$ is defined by the same formula (7.1), and the
right mapping~$r$ is defined by (7.2), 
i.e., 
$r(x,p)=s^{-1}_x(\ell(x,p))$, where $s_x$ is the inversion 
mapping.

\begin{lemma}
The components of the mappings $\ell, r$ defined above 
obey the Lie--Engel system {\rm(7.3)}.
\end{lemma}

The family of inversions $\{s_x\}$ still defines a Lagrangian
fibration of a neighborhood of the diagonal in $\cX\times\cX$. 
The fibers are graphs of the mappings~$s_x$.

If $\gamma$ is a transformation of~$\cX$, 
the $\Graph(\gamma)$ intersects the fiber of $s_x$ 
at a point $(\gamma(\tilde{x}),\tilde{x})$,
where $\tilde{x}$ is a fixed point of the mapping 
$s^{-1}_x\circ\gamma$
(compare with Sec.~3).

\begin{lemma}
Theorem~{\rm3.1,\,(i)} still holds.
If the mapping is symplectic, 
then the form $\cH_x(\gamma(\tilde{x}))$ is closed, 
and we can define a function $\Phi^\gamma$ by the formula
\begin{equation}
D\Phi^\gamma(x)=\cH_x(\gamma(\tilde{x}))
\tag{10.2}
\end{equation}
{\rm(}compare with {\rm(3.2))}.
\end{lemma}

Let us now define 
{\it internal exponential mappings} $\Exp_x$ 
in the same way as in Sec.~2.

Consider the curve $\sigma_x=\sigma^+_x\cup\sigma^-_x$ 
composed of two pieces: 
$\sigma^+_x=\{\Exp_x(vt)\mid 0\leq t\leq 1\}$ and 
$\sigma^-_x=\{s^{-1}_x(\Exp_x(vt))\}$.
Denote $z=\Exp_x(v)$, $y=s^{-1}_x(z)$.
We call $\sigma_x$ the 
{\it internal geodesic form~$y$ to~$z$ trough the center~$x$}.

Since $s_x(\sigma^+_x)=\sigma^-_x$, 
such an internal geodesic
$\sigma_x$ is an inversive curve with the center point~$x$.

\begin{lemma}
{\rm(i)} The membrane formula {\rm(3.3)}
for the phase function of the symplectic transformation~$\gamma$
still holds if in the construction of the membrane 
$\Sigma^\gamma(x,y)$
the Ether geodesics are replaced by the internal geodesics 
and
the term ``mid-point'' is replaced by the term ``center-point.''

{\rm(ii)} The membrane formulas {\rm(4.1), (6.1)},
the Hamilton--Jacobi equation {\rm(7.5)},
the identity {\rm(4.2)}, 
as well as the composition formulas {\rm(6.3), (6.8)}
and all other formulas of Secs.~{\rm7--9},
hold with the same replacement agreement as in case~{\rm(i)}.
\end{lemma}

\bibliographystyle{amsalpha}

\end{document}